\documentclass[12pt]{amsart}
\usepackage{amsmath,amsfonts,amsthm,amscd,}
\newtheorem{theorem}{Theorem}

\newtheorem{definition}[theorem]{Definition}

\newtheorem{corollary}[theorem]{Corollary}
\theoremstyle{remark}

\newcommand{\ds}{\displaystyle}
\newcommand{\sk}{\sigma_k^{1/k}(S^n)}

\numberwithin{equation}{section} \numberwithin{theorem}{section}
\begin{document}
\bibliographystyle{amsalpha}
\title[Conformal geometry]{Conformal geometry and fully nonlinear equations}
\author{Jeff Viaclovsky}
\address{Jeff Viaclovsky, Department of Mathematics, MIT, 
Cambridge, MA 02139}
\address{Department of Mathematics, University of Wisconsin, 
Madison, WI, 53706}
\email{jeffv@math.wisc.edu}
\thanks{The research of the author was partially
supported by NSF Grant DMS-0503506.}
\begin{abstract}
This article is a survey of results involving conformal 
deformation of Riemannian metrics and fully nonlinear equations. 
\end{abstract}
\date{August 31, 2006}
\dedicatory{To the memory of Professor S.S. Chern}
\maketitle
\section{The Yamabe equation}
One of the most important problems in conformal 
geometry is the Yamabe Problem, which 
is to determine whether there exists a conformal metric with 
constant scalar curvature on any closed Riemannian manifold.
In what follows, let $(M,g)$ be a Riemannian manifold, and 
let $R$ denote the scalar curvature of $g$. 
Writing a conformal metric as $\tilde{g} = v^{\frac{4}{n-2}} g$, 
the Yamabe equation takes the form
\begin{align}
\label{Yam}
 4 \frac{n-1}{n-2} \Delta v + R \cdot v = \lambda \cdot v^{\frac{n+2}{n-2}},
\end{align}
where $\lambda$ is a constant. These are the Euler-Lagrange equations 
of the {\em Yamabe functional}, 
\begin{align}
\mathcal{Y}(\tilde{g}) =  Vol(\tilde{g})^{- \frac{n-2}{n}} 
\int_M R_{\tilde{g}} dvol_{\tilde{g}},
\end{align}
for $\tilde{g} \in [g]$, where $[g]$ denotes the conformal class of $g$. 
An important related conformal invariant is the 
{\em Yamabe invariant} of the conformal class $[g]$: 
\begin{align}
Y([g]) \equiv \underset{ \tilde{g} \in [g] }{\mbox{inf }} \mathcal{Y}(\tilde{g}).
\end{align}
 The Yamabe problem has been completely solved through the results 
of many mathematicians, over a period of approximately thirty years. 
Initially, Yamabe claimed to 
have a proof in \cite{Yamabe}.
 The basic strategy was to prove the existence of a minimizer
of the Yamabe functional through a sub-critical regularization technique. 
 Subsequently, an error was
found by N. Trudinger, who then gave a solution with a 
smallness assumption on the Yamabe invariant \cite{TrudingerYam}. 
Later, Aubin showed that
the problem is solvable provided that 
\begin{align}
Y([g]) < Y([g_{round}]),
\end{align}
where $[g_{round}]$ denotes the conformal class of the round metric on the 
$n$-sphere, 
and verified this inequality for $n \geq 6$ and $g$ not locally conformally 
flat \cite{Aubin1}, \cite{Aubin2}, \cite{Aubin3}. Schoen 
solved the remaining cases \cite{Schoen}. 
It is remarkable that Schoen employed the positive mass
conjecture from general relativity to solve these 
remaining most difficult cases. 

An important fact is that $SO(n+1,1)$, 
the group of conformal transformations of 
the $n$-sphere $S^n$ with the round metric, is non-compact.
Likewise, the space of solutions 
to the Yamabe equation in the conformal class of the round sphere 
is non-compact. However, if $(M,g)$ is compact, and not conformally equivalent 
to the round sphere, then the group of conformal 
transformations is compact \cite{Ferrand}, \cite{Obata}. 
A natural question is then whether the space of 
all unit-volume solutions (not just minimizers) to (\ref{Yam}) is compact
on an arbitrary compact manifold, 
provided $(M,g)$ is not conformally equivalent to $S^n$ 
with the round metric. Schoen solved this in the locally conformally
flat case \cite{Schoen3}, and produced unpublished 
lecture notes outlining a solution is certain other cases. 
Many other partial solutions have appeared, see 
for example \cite{Druet}, \cite{LiZhangI}, \cite{LiZhangII}, 
\cite{Marques}, \cite{Schoen1}.
Schoen has recently announced the complete solution of the 
compactness problem in joint work with Khuri and Marques, assuming
that the positive mass theorem holds in higher dimensions.
The positive mass theorem is known to hold in the locally conformally flat case
in all dimensions \cite{SYI}, and in the general case in 
dimensions $n \leq 7$ \cite{SYII}, \cite{LeeParker},
and in any dimension if the manifold is spin \cite{Witten},
\cite{ParkerTaubes}, \cite{LeeParker}. 

\section{A fully nonlinear Yamabe problem}
   The equation (\ref{Yam}) is a semi-linear equation, meaning 
the the non-linearities only appear in lower order 
terms -- second derivatives appear in a linear fashion. 
One may investigate other types of conformal curvature 
equations, which brings one into the realm of fully 
nonlinear equations. We recall the Schouten tensor
\begin{align}
A_g = \frac{1}{n-2} \left(  Ric - \frac{R}{2(n-1)} g \right),
\end{align}
where $Ric$ denotes the Ricci tensor. 
This tensor arises naturally in the decomposition of 
the full curvature tensor
\begin{align}
Riem = Weyl + A \odot g, 
\end{align}
where $\odot$ denotes the Kulkari-Nomizu product 
\cite{Besse}. This equation also serves to define the Weyl
tensor, which is conformally invariant. 
Thus the behaviour of the full curvature tensor 
under a conformal change of metric is entirely 
determined by the Schouten tensor. 
Let $F$ denote any symmetric function of the 
eigenvalues, which is homogeneous of degree one, and consider the equation
\begin{align}
\label{Fyam}
F(\tilde{g}^{-1} A_{\tilde{g}}) = \mbox{constant}.
\end{align}
Note that the $\tilde{g}^{-1}$ factor is 
present since only the eigenvalues of an endomorphism 
are well-defined.  
If we write a conformal metric as $\tilde{g} = e^{-2u} g$, 
the Schouten tensor transforms as 
\begin{align}
A_{\tilde{g}} = \nabla^2 u + du \otimes du - \frac{|\nabla u|^2}{2} g + A_g. 
\end{align}
Therefore, equation (\ref{Fyam}) is equivalent to 
\begin{align}
F \left( g^{-1} \left(
\nabla^2 u + du \otimes du - \frac{| \nabla u|^2}{2} g + A_g
\right) \right ) = \mbox{constant} \cdot e^{-2u}.
\end{align}
Let $\sigma_k$ denote the $k$th elementary symmetric 
function of the eigenvalues
\begin{align}
\sigma_k = \sum_{i_1 < \dots < i_k} \lambda_{i_1} \cdots \lambda_{i_k}.
\end{align}
For the case of $F = \sigma_k^{1/k}$, the equation (\ref{Fyam}) 
has become known as the {\em $\sigma_k$-Yamabe 
equation}:
\begin{align}
\label{skyam}
\sigma_k^{1/k}(\tilde{g}^{-1} A_{\tilde{g}}) = \mbox{constant}.
\end{align}
In the context of exterior differential systems, they 
arose in a different form in Bryant and Griffiths' research 
on conformally invariant Poincar\'e-Cartan forms \cite{BGG}, 
and these systems were shown to correspond to the 
$\sigma_k$-Yamabe equation in \cite{Jeff1}. 

For $1 \leq k \leq n$, we define the cone (in $\mathbb{R}^n$)
\begin{align}
\Gamma_k^+ = \{ \sigma_k > 0 \} \cap \{ \sigma_{k-1} > 0 \}
\cap \cdots \cap \{ \sigma_1 > 0 \}.
\end{align}
These are well-known as ellipticity cones for the 
$\sigma_k$ equation, see \cite{Garding}, 
\cite{Ivochkina}, \cite{CNSIII}. 
We will say that a metric $g$ is strictly $k$-{\em{admissible}}
if the eigenvalues of $g^{-1}A_g$ lie in $\Gamma_k^+$ 
at every point $p \in M$. It is an important fact 
that if the metric $\tilde{g}$ is $k$-admissible, then the linearization of 
(\ref{skyam}) at $\tilde{g}$ is elliptic. 
On a compact manifold, if the background metric $g$ is 
$k$-admissible, then (\ref{skyam}) is necessarily elliptic at {\em{any}} 
solution \cite{Jeff2}. Thus, a $k$-admissiblity 
assumption on the background metric is an ellipticity assumption, 
reminiscent of the $k$-convexity assumption on domains for the
$k$-Hessian equation \cite{CNSIII}. 
\section{Variational characterization}
In \cite{Jeff1}, it was shown that 
in several cases, the $\sigma_k$-Yamabe equation is 
variational. Let $\mathcal{M}_1$ denote the set of unit volume
metrics in the conformal class $[g_0]$.
 \begin{theorem}(Viaclovsky \cite{Jeff1})
\label{ELEquations}
If $k \neq n/2$ and $(N , [g_0])$ is locally conformally flat,
a metric $g \in \mathcal{M}_1$ is a critical point of the functional
$$\mathcal{F}_k : g \mapsto
\int_N \sigma_k ( g^{-1}A_g ) \mbox{dV}_g$$
restricted to $\mathcal{M}_1$
if and only if
$$ \sigma_k ( g^{-1} A_g ) = C_k$$
for some constant $C_k$.
If $N$ is not locally conformally flat, the statement is
true for $k=1$ and $k=2$.
\end{theorem}

For $k=1$, this is of course well-known, as 
$\mathcal{F}_1$ is the Hilbert functional, \cite{Hilbert}, \cite{Schoen1}. 
For $k = n/2$, in \cite{Jeff1} it was shown that the integrand 
is the non-Weyl part of the Chern-Gauss-Bonnet integrand. 
Therefore, $\mathcal{F}_{n/2}$ is necessarily constant in 
the locally conformally flat case 
(when $n=4$, this holds in general, this will 
be discussed in detail in Section \ref{d4} below).
Nevertheless, in this case the equation is still 
variational, but with a different functional. 
Fix a background metric $h$, write $g = e^{-2u}h$, and let
\begin{align}
\mathcal{E}_{n/2}(g) = \int_M \int_0^1 \sigma_{n/2} \left(-t \, \nabla_h^2 u + 
t^2 \, \nabla_h u 
\otimes \nabla_h u - \frac{1}{2} \, t^2 \, 
|\nabla_h u|^2 \, g_0 + A(h) \right) \, u \, dt \, dV_h,
\end{align}
then for any differentiable path of smooth conformal metrics $g_t$, 
\begin{align}
\frac{d}{dt} \mathcal{E}_{n/2}(g_t) = \int_M \sigma_{n/2}(A_{g_t}) \, u \,
dvol_{g_t}.
\end{align}
This fact was demonstrated in \cite{BrendleViaclovsky},
see also \cite{ChangYang}. This is valid also for 
$n=4$, \cite{ChangYangAnnals}.
Recently, Sheng-Trudinger-Wang have given 
conditions on when the more general $F$-Yamabe equation 
is variational, see \cite{STW}.

A natural question is: what are the 
critical metrics of the $\mathcal{F}_k$ 
functionals, when considering all possible 
metric variations, not just conformal variations?
It is a well-known result that the critical points of
$\mathcal{F}_1$ restricted to space of unit volume 
metrics are
exactly the Einstein metrics \cite{Schoen1}. But for $k > 1$, the 
full Euler-Lagrange equations are manifestly 
fourth order equations in the metric. 
However, in dimension three we have the following
\begin{theorem}\label{maintheorem} (Gursky-Viaclovsky \cite{GurskyViaclovsky}) 
Let $M$ be compact and of dimension three. Then a metric $g$ with 
$\mathcal{F}_2[g] \geq 0$ is critical for 
$\mathcal{F}_2$ restricted to the space of 
unit volume metrics if and only if $g$ has constant sectional curvature. 
\end{theorem}
 A similar theorem was proved by Hu-Li for $ n \geq 5$, but with the 
rather stringent condition that the metric be 
locally conformally flat \cite{HuLi}.
We mention that Labbi studied some
curvature quantities defined by H. Weyl which are polynomial in the 
full curvature tensor, and proved some 
interesting variational formulas \cite{Labbi}. 
\section{Liouville Theorems}
We next turn to the uniqueness question. In the negative 
curvature case, the linearization of (\ref{Fyam}) is 
invertible, so the uniqueness question is trivial. 
However, in the positive curvature case the uniqueness 
question is non-trivial.
In \cite{Jeff1}, \cite{JeffT} the following was
proved
\begin{theorem}(Viaclovsky \cite{Jeff1})
\label{uniqueness}
Suppose $(N,g_0)$ is of unit volume and has constant sectional
curvature $K > 0$.
Then for any $k \in \{1, \dots, n \}$, $g_0$ is the
unique unit volume solution in its conformal class of
\begin{align}
\label{Ricci}
\sigma_k ( g^{-1} A_g ) = {\mbox{constant}},
\end{align}
unless $N$ is isometric to $S^n$ with the round metric.
In this case we have an $(n+1)$-parameter family of solutions
which are the images of the standard metric under the conformal
diffeomorphisms of $S^n$.
\end{theorem}
For $k=1$, the constant scalar curvature case, the
theorem holds just assuming $N$ is Einstein.
This is the well-known theorem of Obata \cite{Obata}.

This theorem falls under the category of a Liouville-type 
theorem.
We let $\delta_{ij}$ be the Kronecker delta symbol, 
and write the conformal factor as $\tilde{g} = u^{-2} g$. 
In stereographic coordinates, the equation (\ref{Ricci}) is written 
\begin{align}
\label{matrix}
\sigma_k\left( u \cdot \frac{\partial^2u}{\partial{x^i}\partial{x^j}} - \frac{|\nabla u|^2}{2}\delta_{ij}
\right) = {\mbox{constant}}.
\end{align}
This equation is conformally invariant: 
if $T : \mathbb{R}^n \rightarrow \mathbb{R}^n$ is
a conformal transformation (i.e., $T \in SO(n+1,1)$), 
and $u(x)$ is a solution of (\ref{matrix}),
then 
\begin{align}
\label{invariance}
v(x) = |J(x)|^{-1/n} u(Tx)
\end{align}
is also a solution, where $J$ is the Jacobian 
of $T$, see \cite{JeffT}. 

 The uniqueness theorem can then be restated as a Liouville Theorem
in $\mathbb{R}^n$:
\begin{theorem}({Viaclovsky} \cite{JeffT})
Let $u(x) \in C^{\infty}( \mathbb{R}^n)$ be a positive solution to
\begin{align}
\label{matrix2}
\sigma_k\left( u \cdot \frac{\partial^2u}{\partial{x^i}\partial{x^j}} - \frac{|\nabla u|^2}{2}\delta_{ij}
\right) = {\mbox{constant}}
\end{align}
for some $k \in \{1, \dots, n \}$.
Suppose that $v(y) = |y|^2 \cdot u(\frac{y^1}{|y|^2}, \dots,
\frac{y^n}{|y|^2})$ is smooth and
$$ \lim_{y \rightarrow 0} v(y) > 0.$$ Then
$$u = a|x|^2 + b_ix^i + c,$$
where $a$, $b_i$, and $c$ are constants.
\end{theorem}
The proof in \cite{Jeff1}, \cite{JeffT} requires the stringent growth 
condition at infinity. 
For the scalar curvature equation, $k=1$, this 
theorem was proved without {\em any} assumption 
at infinity in the important paper by 
Caffarelli-Gidas-Spruck \cite{CGS},
using the moving planes technique. 

 This analagous theorem for $ k \geq 2$ is 
now known to hold without any condition on the 
behaviour at infinity. 
For $k=2$, important work was done in \cite{CGYA}, 
and \cite{CGYE}, proving the Liouville Theorem 
for $n = 4,5$, and for $ n \geq 6$ with a finite 
volume assumption. Maria del Mar Gonz\'alez proved 
a Liouville Theorem for $\sigma_k$, $n > 2(k+1)$ 
with a finite volume assumption \cite{Gonzalez2}.
 
Yanyan Li and Aobing Li proved the following theorem for 
all $k$ in the important paper \cite{LiLi}.
\begin{theorem}(Li-Li \cite{LiLi})
Let $u(x) \in C^{\infty}( \mathbb{R}^n)$ be a positive solution to
\begin{align}
\label{matrix3}
\sigma_k\left( u \cdot \frac{\partial^2u}{\partial{x^i}\partial{x^j}} - \frac{|\nabla u|^2}{2}\delta_{ij}
\right) = {\mbox{constant}} > 0
\end{align}
for some $k \in \{1, \dots, n \}$, 
satisfying 
\begin{align}
\label{conecond}
u \cdot \frac{\partial^2u}{\partial{x^i}\partial{x^j}} - \frac{|\nabla u|^2}{2}\delta_{ij}\in \Gamma_k^+.
\end{align}
Then
$$u = a|x|^2 + b_ix^i + c,$$
where $a$, $b_i$, and $c$ are constants.
\end{theorem}
Their method is based on the moving planes technique. 
Subsequently, their results have been 
generalized to much more general classes of 
symmetric functions $F$, see \cite{LiLiII}, 
\cite{Li3}, \cite{LiLi4}, ,\cite{Li5}.

 These types of Liouville theorems can be used in deriving 
a priori estimates for solutions of (\ref{Ricci}), we 
will discuss this below.

\section{Local estimates}
\label{lest}
We consider the $\sigma_k$-Yamabe equation 
\begin{align}
\label{kyam}
\sigma_k^{1/k} \left( \nabla^2 u + du \otimes du - \frac{ |\nabla u|^2}{2} g
+ A_g \right) = f(x) e^{-2u},
\end{align}
with $f(x) \geq 0$. A remarkable property of (\ref{kyam}) 
was discovered in \cite{GuanWang1}. It turns out that
{\em local } estimates are satisfied, a fact which does 
not hold in general fully nonlinear equations. 
We say that $u \in C^2$ is $k-$admissible if $A_u \in
\overline{\Gamma}_{k}^+$.
For equation (\ref{kyam}), Guan and Wang prove
\begin{theorem}(Guan-Wang \cite{GuanWang1})
\label{GWlocal}
Let $u \in C^3(M^n)$ be a $k-$admissible solution of (\ref{kyam})
in $B(x_0,\rho)$, where $x_0 \in M^n$ and $\rho
> 0$.  Then there is a constant
\begin{align*}
C_0 = C_0(k,n,\rho,\|g\|_{C^2(B(x_0,\rho))},
\|f\|_{C^1(B(x_0,\rho))}),
\end{align*}
such that
\begin{align} \label{GWlocalest}
|\nabla u|^2(x) \leq C_0\big(1 + e^{-2 \inf_{B(x_0,\rho)}u}\big)
\end{align}
for all $x \in B(x_0,\rho/2)$.

Let $u \in C^4(M^n)$ be a $k-$admissible solution of (\ref{kyam})
in $B(x_0,\rho)$, where $x_0 \in M^n$ and $\rho
> 0$.  Then there is a constant
\begin{align*}
C_0 = C_0(k,n,\rho,\|g\|_{C^3(B(x_0,\rho))},
\|f\|_{C^2(B(x_0,\rho))}),
\end{align*}
such that
\begin{align} \label{GWlocalest2}
|\nabla^2 u|(x) + |\nabla u|^2(x) \leq C_0\big(1 + e^{-2
\inf_{B(x_0,\rho)}u}\big)
\end{align}
for all $x \in B(x_0,\rho/2)$.

\end{theorem}

These local estimates for (\ref{kyam}) generalize the global estimates 
which were first proved in \cite{Jeff2}.
Subsequently, these results have been extended to much more general 
classes of symmetric functions $F$, see \cite{Sophie1}, \cite{GuanWangDuke}, 
\cite{GuanLinWang2}, \cite{LiLi},\cite{Yanyanlocal}, \cite{WangXJ}.
Estimates for solutions of $\sigma_2$ in dimension four were proved in 
\cite{Hanlocal} using integral methods.

  Equipped with second derivative estimates, one then uses the 
work of Evans and Krylov \cite{Evans}, \cite{Krylov}
to obtain $C^{2, \alpha}$ estimates, that is, a H\"older estimate 
on second derivatives. This is crucial -- the importance
of that work in this theory can not be overstated.  

Consider a symmetric function
\begin{align} \label{Fex}
F : \Gamma \subset \mathbb{R}^n \rightarrow \mathbb{R}
\end{align}
with $F \in C^{\infty}(\Gamma) \cap C^{0}(\overline{\Gamma})$, where
$\Gamma \subset \mathbb{R}^n$ is an open, symmetric, convex cone,
and impose the following conditions:

\vskip.1in $(i)$ $F$ is symmetric, concave, and homogenous of degree
one.

\vskip.1in  $(ii)$ $F > 0$ in $\Gamma$, and $F = 0$ on $\partial
\Gamma$.

\vskip.1in  $(iii)$  $F$ is {\em elliptic}: $F_{\lambda_i}(\lambda)
> 0$ for each $1 \leq i \leq n$, $\lambda \in \Gamma$.

\vskip.1in

We mention that Szu-Yu Chen proved local $C^2$-estimates (\ref{GWlocalest2}) 
for this general class of symmetric functions $F$ \cite{Sophie1}.

An immediate corollary of these local estimates is an $\epsilon-$regularity
result:

\begin{theorem}(Guan-Wang \cite{GuanWang1})
\label{epsilonregr} There exist constants $\epsilon_0 > 0$ and $C =
C(g,\epsilon_0)$ such that any solution $u \in C^2(B(x_0,\rho))$ of
(\ref{kyam}) with
\begin{align}  \label{smallrenergy}
\int_{B(x_0,\rho)} e^{-nu} dvol_g \leq \epsilon_0,
\end{align}
satisfies
\begin{align} \label{lowerb}
\inf_{B(x_0,\rho/2)} u \geq -C + \log \rho.
\end{align}  Consequently, there is a constant
\begin{align*}
C_2 = C_2(k,n,\mu,\epsilon_0,\|g\|_{C^3(B(x_0,\rho))}),
\end{align*}
such that
\begin{align} \label{gradonr}
|\nabla^2 u|(x) + |\nabla u|^2(x) \leq C_2\rho^{-2}
\end{align}
for all $x \in B(x_0,\rho/4)$.
\end{theorem}

This type of estimate is crucial in understanding 
{\em bubbling}, a phenomenon which is unavoidable when 
studying conformally invariant problems. It
shows that non-compactness of the space 
of solutions can arise only through volume concentration. 
\section{Dimension four}
\label{d4}
In dimension four, an important conformal invariant is 
\begin{align}
\label{F2}
\mathcal{F}_2([g])\equiv 4 \int_M \sigma_2(A_g) dV_g =
\int_M \left( -\frac{1}{2} |Ric_g|^2 + \frac{1}{6} R_g^2 \right) dvol_g.
\end{align}  
By the Chern-Gauss-Bonnet formula (\cite{Besse}),
\begin{align}
\label{CGB}
8\pi^2\chi(M) = \int_M |W_g|^2 dvol_g  + \mathcal{F}_2([g]).
\end{align}
Thus, the conformal invariance of $\mathcal{F}_2$ follows from the well known
(pointwise) conformal invariance of the Weyl tensor $W_g$ 
(see \cite{Eisenhart}).  

One of the most interesting results in this area is 
\begin{theorem}(Chang-Gursky-Yang \cite{CGY1})
\label{AMP}
Let $(M,g)$ be a closed 4-dimensional Riemannian 
manifold with positive scalar curvature. If
$\mathcal{F}_2([g]) > 0$, 
then there exists a conformal metric $\tilde{g} = e^{-2u} g$ with
$R_{\tilde{g}} > 0$ and $\sigma_2(\tilde{g}^{-1}A^1_{\tilde{g}}) > 0$ pointwise.  In particular, the Ricci curvature of 
$\tilde{g}$ satisfies 
$$
0 < 2 Ric_{\tilde{g}} <  R_{\tilde{g}} \tilde{g}.
$$
\end{theorem}

 By combining this with some work of Margerin \cite{Margerin} on the 
Ricci flow, the authors obtained the 
following remarkable integral sphere-pinching theorem:
\begin{theorem}(Chang-Gursky-Yang \cite{CGYSphere}
Let $(M^4, g)$ be a smooth, closed four-manifold for which 
\begin{itemize}
\item[{\rm (i)}] the Yamabe invariant $Y([g]) > 0$, and 
\item[{\rm (ii)}] the Weyl curvature satisfies
\begin{align}
\int_{M^4} |W|^2 dvol < 16 \pi^2 \chi(M^4).
\end{align}
\end{itemize}
Then $M^4$ is difeomorphic to either $S^4$ or $\mathbb{RP}^4$. 
\end{theorem}

The proof of Theorem \ref{AMP} in \cite{CGY1} involved regularization by a 
fourth-order equation and relied on some delicate integral estimates.  
Subsequently, a more direct proof was given in \cite{GVJDG}:
define the tensor
\begin{align}
A^t_g = \frac{1}{2} \Big( Ric_g - \frac{t}{6} R_g g \Big).
\end{align}
\begin{theorem}
\label{main}
(Gursky-Viaclovsky \cite{GVJDG}) Let $(M,g)$ be a closed 4-dimensional 
Riemannian manifold with positive scalar curvature. If  
\begin{align} 
\label{assump}
 \mathcal{F}_2([g]) + \frac{1}{6} (1-t_0) (2-t_0) (Y([g]))^2 > 0,
\end{align}
for some $t_0 \leq 1$, 
then there exists a conformal metric $\tilde{g} = e^{-2u} g$ with
$R_{\tilde{g}} > 0$ and $\sigma_2(A^{t_0}_{\tilde{g}}) > 0$ pointwise. 
This implies 
the pointwise inequalities
\begin{align}
\label{Ricineq}
(t_0-1) R_{\tilde{g}} \tilde{g} < 2 Ric_{\tilde{g}} < (2-t_0) 
R_{\tilde{g}} \tilde{g}.
\end{align}
\end{theorem}

The proof involves a deformation of the equation through 
a path of fully nonlinear equations, and an application 
of the local estimates of Guan-Wang. 
A similar technique was applied the $\sigma_k$ equations in the 
locally conformally flat case by Guan-Lin-Wang \cite{GuanLinWang}. 

 As applications of Theorem \ref{main}, consider two different 
values of $t_0$.  When $t_0 = 1$, we obtain 
the aforementioned result in \cite{CGY1}.
The second application is to the spectral properties of a conformally 
invariant differential operator known as the 
{\it Paneitz operator}.  Let $\delta$ denote the $L^2$-adjoint of the 
exterior derivative $d$; then the Paneitz operator is
defined by 
\begin{align}
P_{g} \phi   = \Delta^2 \phi + \delta \Big( \frac{2}{3} R_g g - 2 Ric_g \Big)
d\phi. 
\end{align}
The Paneitz operator is conformally invariant, in the sense
that if $\tilde{g} = e^{-2u}g$, then
\begin{align}
\label{confinv}
P_{\tilde{g}} = e^{4u}P_g.
\end{align}
Since the volume form of the conformal metric $\tilde{g}$ is $dvol_{\tilde{g}} = e^{-4u}dvol_g$, an immediate consequence of 
(\ref{confinv}) is the conformal invariance of the Dirichlet energy 
\begin{align*}
\langle P_{\tilde{g}} \phi, \phi \rangle_{L^2(M,\tilde{g})} = \langle P_g \phi, \phi \rangle_{L^2(M,g)}.
\end{align*}
In particular, positivity of the Paneitz operator is a conformally invariant property, and clearly the kernel is invariant as well.

To appreciate the geometric significance of the Paneitz operator, define the 
associated $Q$-{\it curvature}, introduced by Branson: 
\begin{align} 
Q_g = - \frac{1}{12} \Delta R_g  + 2 \sigma_2(g^{-1}A^1_g).
\end{align}
Under a conformal change of metric 
$\tilde{g} = e^{-2u}g$, the $Q$-curvature transforms according to the 
equation
\begin{align}
\label{Qchange}
-Pu + 2Q_g = 2Q_{\tilde{g}}e^{-4u},
\end{align}
see, for example, \cite{BransonOrsted}.
Note that 
\begin{align}
\label{intQ}
\int_M Q_g dvol_g = \frac{1}{2} \mathcal{F}_2([g]),
\end{align}
so the integral of the $Q$-curvature is conformally invariant.

An application of Theorem \ref{main} with
$t_0=0$, yields
\begin{theorem}(Gursky-Viaclovsky \cite{GVJDG})
\label{Q}
Let $(M,g)$ be a closed 4-dimensional Riemannian 
manifold with positive scalar curvature. If
\begin{align} 
\label{assump3}
 \int Q_g dvol_g + \frac{1}{6}(Y[g])^2 > 0,
\end{align}
then the Paneitz operator is nonnegative, and 
$Ker P = \{ constants \}$. Therefore, 
by the results in \cite{ChangYangAnnals}, there 
exists a conformal metric $\tilde{g} = e^{-2u} g$
with $Q_{\tilde{g}} = constant.$
\end{theorem}
This yields new examples of manifolds admitting 
constant $Q$-curvature metrics, see \cite{GVJDG}. 

 Another interesting application of $\sigma_2$ in dimension
$4$ was found in \cite{WangDirac}, which gives a lower
estimate on eigenvalues of the Dirac operator 
in terms of $\mathcal{F}_2$ on a spin manifold.

\section{Parabolic methods}
We recall the {\em Yamabe flow},
\begin{align}
\frac{d}{dt} g = -(R_g - r_g) g,
\end{align}
where $r_g$ denotes the mean value of the 
scalar curvature. This flow was introduced by Hamilton, 
who proved existence of the flow for all time
and proved convergence in the case of negative scalar curvature. 
The case of positive scalar curvature however is highly 
non-trivial. The locally conformally flat case was studied in 
\cite{ChowJDG} and \cite{YeJDG}. 
Schwetlick and Struwe \cite{SS} proved convergence
for $3 \leq n \leq 5$ provided an certain energy 
bound on the initial metric is satsified. 
In the beautiful paper \cite{Brendle}, Simon Brendle 
proved convergence for $3 \leq n \leq 5$ 
for {\em{any}} initial data. 

In the fully nonlinear case, the following flow was first 
proposed in \cite{GuanWang2}:
\newcommand{\ba}{\begin{array}}
\newcommand{\ea}{\end{array}}
\def \ds{\displaystyle}
\def \vs{\vspace*{0.2cm}} 
\begin{equation}\label{flow1}
\left\{\ba{rcl}
\ds\vs \frac {d}{dt}g
&=&-(\log \sigma_k(g)-\log r_k(g))\cdot g, \\
g(0)&=&g_0,\ea\right.
\end{equation}
where $r_k(g)$ is given by
\[ r_k(g)=  \exp \left(\frac 1{vol(g)} \int_M \log \sigma_k (g)
\, dvol(g)\right).\]

If $g=e^{-2u}\cdot g_0$, then equation (\ref{flow1}) can be written 
as the following fully nonlinear flow
\def\uuu{{\nabla^2 u+du\otimes du-\frac {|\nabla u|^2} 2 g_0+A_{g_0}}}
\begin{equation}\label{flow2}
\left\{\ba{rcl}    
\ds \vs 2\frac{\ds du}{\ds dt} & = &  \ds \log \sigma_k
\left(\uuu\right)+2ku -\log r_k\\
u(0) & = & u_0.
\ea\right.
\end{equation}
Guan and Wang settled the locally conformally flat case:
\begin{theorem}(Guan-Wang \cite{GuanWang2})\label{thm2}
Suppose $(M,g_0)$ be a compact, connected and locally conformally flat
manifold. Assume that $A_{g_0} \in \Gamma_k^+$ and smooth, then flow
(\ref{flow1}) exists for all time $0<t <\infty$ and 
$g(t) \in C^{\infty}(M)$ for all $t$.
For any positive integer $l$, there exists a constant $C$ 
depending only on $g_0, k, n$ (independent of $t$) such that
\begin{align}
\Vert g \Vert_{C^l(M)} \leq C,
\end{align}
where the norm is taken with respect to the background 
metric $g_0$. Furthermore, there are a positive number $\beta$ and a smooth 
metric $g_{\infty} \in \Gamma_k^+$ such that
\begin{align}
\sigma_k (A_{g_{\infty}}) = \beta, 
\end{align}
and
\begin{align}
\lim_{t \rightarrow \infty} \Vert g(t) - g_{\infty} \Vert_{C^1(M)} = 0,
\end{align}
for all $l$. 
\end{theorem}
Guan and Wang employed a log flow (rather than a gradient 
flow) due to a technical reason in obtaining $C^2$ estimates. 
This solved the $\sigma_k$-Yamabe problem in the 
locally conformally flat case. Independently, Yanyan Li
and Aobing Li solved the locally conformally flat 
$\sigma_k$-Yamabe problem using elliptic methods \cite{LiLi}. 

In a subsequent paper \cite{GuanWangDuke}, Guan-Wang employed 
the log-flow for the quotient equations to obtain some 
interesting inequalities.
Define the scale invariant functionals
\begin{align}
\mathcal{F}_k(g) = (Vol(g))^{- \frac{n - 2k}{n}}
\int_M \sigma_k ( g^{-1} A_g) dV_g. 
\end{align}
\def\cal{\mathcal}
\def\Fk{{\cal F}_k}
\def\s{\sigma}
\def\sk{\s_k}
\def\Fl{{\cal F}_l}
\def\lcf{{locally conformally flat} }
\begin{theorem}(Guan-Wang \cite{GuanWangDuke})
\label{inqs} Suppose that $(M,g_0)$ is a compact, oriented and connected
\lcf manifold with
$A_{g_0} \in \Gamma_k^+$ smooth and conformal metric $g$ 
with $A_g \in \Gamma_k^+$. Let $0\le l< k\le n$.
\begin{itemize}
\item[({A}).] {\em Sobolev type inequality:} If
$0\le l < k < \frac n2$, then there is a positive constant $C_S=C_S([g_0], n,k,l)$
depending only on  $n$, $k$, $l$ and the
conformal class $[g_0]$ such that
\begin{equation}\label{Sob-ineq} \left(\Fk(g)\right)^{\frac 1{n-2k}}
\ge C_S \left(\Fl(g)\right)^{\frac 1{n-2l}}.\end{equation} 
\item[({B}).] {\em Conformal quermassintegral type
inequality:} If $n/2\le k \le n$, $1\le l<k$, then
\begin{equation}\label{cquer}
\left(\Fk(g)\right)^{\frac 1 k}\le
{\binom{n}{k}}^{\frac 1k}{\binom{n}{l}}^{-\frac 1l} \left(\Fl(g)\right)^{\frac 1l}.\end{equation}
\item[({C}).] {\em Moser-Trudinger type inequality:}
If $k=n/2$, then
\begin{equation}\label{MT-ineq}
(n-2l){\cal E}_{n/2}(g) \ge C_{MT} \left\{
\log\int_M \s_l(g^{-1}A_g)dV_g-\log 
\int_M \s_l(g_0^{-1} A_{g_0})dV_{g_0}\right\},\end{equation}
\end{itemize}
\end{theorem}

 Furthermore, the constants are all explicit and optimal, with 
a complete characterization of the case of equality see \cite{GuanWangDuke}. 

In the non-locally conformally flat case, recall from Theorem \ref{ELEquations}
that the $\sigma_2$ equation is always variational. 
Using this fact, the $\sigma_2$-Yamabe equation has 
recently been studied in the general case using parabolic 
methods. 
In \cite{GeWang}, convergence 
was proved in dimensions $n > 8$, by constructing an 
explicit test function. In \cite{STW}, 
convergence was proved in dimensions $n > 4$, who 
avoided having to directly construct a test function 
by employing the solution of the Yamabe problem. 
Combining this with the work described in 
Section \ref{posRicsec}, it follows that the $\sigma_2$-Yamabe problem has 
been solved in all dimensions. 
Subsequently, the quotient equation $\sigma_2 / \sigma_1$ 
was studied in \cite{GeWang2} and existence was proved in 
dimensions $n > 4$. 

\section{Positive Ricci curvature}
\label{posRicsec}
The ellipticity assumption of $k$-admissibility 
has geometric consequences on the Ricci curvature.
The following inequality was demonstrated in \cite{GVW}: 
\begin{theorem}(Guan-Viaclovsky-Wang \cite{GVW}) 
\label{kbig} Let $(M,g)$ be a Riemannian manifold
and $x \in M$. If $A_g \in \Gamma_k^+$ at $x$ for some $k\ge
 n/2$, then its Ricci curvature is positive at $x$.
Moreover,  if $A_g \in \overline{\Gamma}_k^+$ for some  $k>1$, then
\[ Ric_g \ge \frac{2k-n}{2n(k-1)} R_g \cdot g .\]
In particular if $k\ge \frac n2$, 
\[ Ric_g \ge \frac{(2k-n)(n-1)}{(k-1)} 
{\binom{n}{k}}^{-\frac 1k} 
\sigma^{\frac 1k}_k(A_g) \cdot g .\]
\end{theorem}
Therefore, in the case $k > n/2$, any strictly $k$-admissible metric 
necessarily has positive Ricci curvature. 
This fact led to  the following definition in 
\cite{GVAM}:
\begin{definition}
Let $(M^n,g)$ be a compact $n$-dimensional Riemannian manifold. For $n/2 \leq k \leq n$ 
the {\em $k$-maximal volume} of $[g]$ is 
\begin{align}
\label{Lambdak}
\Lambda_k(M^n,[g]) = \sup \{vol(e^{-2u}g) | e^{-2u}g \in \Gamma_k^{+}(M^n)\mbox{ with }
\sigma_k^{1/k}(g_u^{-1}A_u) \geq \sk \}.
\end{align}
If $[g]$ does not admit a $k$-admissible metric, set $\Lambda_k(M^n,[g]) = +\infty$.
\end{definition}
Consider the $\sigma_k$-Yamabe equation 
\begin{align}
\label{skYamabe}
\sigma_k^{1/k} \left( \nabla^2 u + du \otimes du - \frac{ |\nabla u|^2}{2} g
+ A_g \right) = C_{n,k} e^{-2u}.
\end{align} 
Where $C_{n,k}$ is the corresponding $\sigma_k^{1/k}$-curvature 
of the standard round metric on $S^n$. 
Using Bishop's inequality, it follows that the invariant
$\Lambda_k$ is non-trivial when $k > n/2$,
and in analogy with the classical Yamabe problem, when the invariant is 
strictly less than 
the value
obtained by the round metric on the sphere one obtains existence of solutions 
to (\ref{skYamabe}):
\begin{theorem}(Gursky-Viaclovsky \cite{GVAM})
\label{Lambdafinite}
If $[g]$ admits a $k$-admissible metric with $k > n/2$, then there is a constant $C = C(n)$ such that
$\Lambda_k(M^n,[g]) < C(n)$.
If 
\begin{align}
\label{subcrit}
\Lambda_k(M^n,[g]) < vol(S^n),
\end{align}
where $vol(S^n)$ denotes the volume of the round sphere,  
then $[g]$ admits a solution $g_u = e^{-2u}g$ of (\ref{skYamabe}).  
Furthermore, the set of solutions of (\ref{skYamabe}) is compact in the 
$C^{m}$-topology for any $m \geq 0$.
\end{theorem}
 The compactness proof uses a bubbling argument, combined 
with the Liouville Theorem of Li-Li \cite{LiLi}, and 
existence is obtained using a degree-theoretic argument. 

In dimension three, we have the following estimate
\begin{theorem}(Gursky-Viaclovsky \cite{GVAM}) 
\label{3dpi1}
Let $(M^3,g)$ be a closed Riemannian three-manifold, and assume $[g]$ admits a
$k$-admissible metric with $k = 2$ or $3$.  Let $\pi_1(M^3)$ denote the fundamental
group of $M^3$.  Then
\begin{align}
\label{sharp3wpi1}
\Lambda_k(M^3,[g]) \leq \frac{vol(S^3)}{\|\pi_1(M^3)\|}.
\end{align}
\end{theorem}
The proof of this used an improvement of Bishop's volume 
comparison theorem in dimension three, due to 
Hugh Bray \cite{Bray}. Therefore, if the three-manifold
is not simply-connected, the estimate (\ref{subcrit}) is 
automatically satisfied. 

In dimension four, the optimal estimate of $\Lambda_k$
follows from the sharp integral estimate for $\sigma_2(A)$
due to Gursky \cite{Gursky1}:

\begin{theorem}(Gursky-Viaclovsky \cite{GVAM})
\label{4dest}
Let $(M^4,g)$ be a closed Riemannian four-manifold, and assume $[g]$ admits a
$k$-admissible metric with $2 \leq k \leq 4$.  Then
\begin{align}
\label{sharp4}
\Lambda_k(M^4,[g]) \leq vol(S^4).
\end{align}
Furthermore, equality holds in (\ref{sharp4}) if and only if $(M^4,g)$ is conformally
equivalent to the round sphere.
\end{theorem}

When $k=2$ the existence was established previously in \cite{CGYA}.
Combining this work with the four-dimensional
solution of the Yamabe problem \cite{Schoen}, it follows that the
$\sigma_k$-Yamabe problem is completely solved in dimension four.

 We return to the case of general dimension $n$. 
Using similar techniques as in the proof of 
Theorem \ref{Lambdafinite}, Guan-Wang studied the minimum 
eigenvalue of the Ricci, and the $p$-Weitzenbock operator, 
and proved various existence theorems, see \cite{GuanWangRic}.
Their problem encountered some additional technical
difficulties since their symmetric function $F$ is not smooth, 
and also required an application of some work of 
Caffarelli \cite{Caffint}, and some generalizations 
of this work \cite{WangLiheI}, \cite{WangLiheII}. 

 In recent work, an existence theorem for a much more 
general class of $F$ was proved. In addition 
to the structure conditions  $(i)-(iii)$ imposed on $F$ above 
in Section \ref{lest}, suppose that we also have
\vskip.1in  $(iv)$ $\Gamma \supset \Gamma_n^{+}$, and there exists a
constant $\delta> 0$ such that any $\lambda =
(\lambda_1,\dots,\lambda_n) \in \Gamma$ satisfies
\begin{align} \label{Ricposcond}
\lambda_i > -\frac{(1- 2\delta)}{(n-2)}\big(\lambda_1 + \dots +
\lambda_n \big) \quad \forall \ 1 \leq i \leq n.
\end{align}

\vskip.1in

To explain the significance of (\ref{Ricposcond}), suppose the
eigenvalues of the Schouten tensor $A_g$ are in $\Gamma$ at each
point of $M^n$.  Then $(M^n,g)$ has positive Ricci curvature: in
fact,
\begin{align} \label{Ricposintro}
Ric_g - 2\delta \sigma_1(A_g)g \geq 0.
\end{align}
Define
\begin{align}
A_u \equiv \nabla^2 u + du \otimes du - \frac{ |\nabla u|^2}{2} g
+ A_g.
\end{align}
For $F$ satisfying $(i)-(iv)$, consider the equation
\begin{align} 
\label{hessF}
F(g^{-1}A_u) = f(x)e^{-2u},
\end{align}
where we assume $g^{-1}A_u \in \Gamma$ (i.e., $u$ is {\em
$\Gamma$-admissible}). 

\begin{theorem}(Gursky-Viaclovsky \cite{GVRic})  \label{Main2}
Suppose $F: \Gamma \rightarrow \mathbb{R}$ satisfies $(i)-(iv)$. Let
$(M^n,g)$ be closed $n$-dimensional Riemannian manifold, and assume
\vskip.1in \noindent $(i)$ $g$ is $\Gamma$-admissible, and
\vskip.1in \noindent $(ii)$ $(M^n,g)$ is not conformally equivalent
to the round $n$-dimensional sphere. \vskip.1in Then given any
smooth positive function $f \in C^{\infty}(M^n)$ there exists a
solution $u \in C^{\infty}(M^n)$ of
\begin{align*}
F(g^{-1}A_u) = f(x)e^{-2u},
\end{align*}
and the set of all such solutions is compact in the $C^m$-topology
for any $m \geq 0$.
\end{theorem}

 This theorem in particular completely solved the $\sigma_k$-Yamabe problem 
whenever $k > n/2$. The proof of this theorem involved an bubbling analysis, 
and the application of various Holder and integral estimates to determine the 
growth rate of solutions at isolated singular points.
For other works analyzing the possible behaviour of solutions 
at singularities, see \cite{CHY}, \cite{GVIMRN}, \cite{Gonzalez2},\cite{Gonzalez3},
\cite{Yanyana}, \cite{Yanyanb}, \cite{TWnew}. 
In addition, a remarkable Harnack inequality for $k$--admissible metrics,
$k > n/2$, was demonstrated in \cite{TWnew}.  

 Note that the second assumption in the above theorem 
is of course necessary, since the set
of solutions of (\ref{hessF}) on the round sphere with $f(x) =
constant$ is non-compact, while for variable $f$ there are
obstructions to existence.  In particular, there is a ``Pohozaev
identity'' for solutions of (\ref{hessF}) in the case
of $\sigma_k$, which holds in the
conformally flat case; see \cite{Jeff3}. 
\begin{theorem}(Viaclovsky \cite{Jeff3})
Let $(M,g)$ be a closed locally conformally flat $n$-dimensional
manifold. Then
for any conformal Killing vector field $X$, 
and $1 \leq k < n$, we have
\begin{align}
\int_M X \cdot \sigma_k(A) { \ dvol_M}
= 0.
\end{align}
\end{theorem}
For $k=1$, this identity is well-known and holds without
the locally conformally flat assumption \cite{SchoenCPAM}.
This identity yields
non-trivial Kazdan-Warner-type obstructions to existence (see
\cite{KazdanWarner}) in the case $(M^n,g)$ is conformally equivalent
to $(S^n, g_{round})$. We note that similar Pohozaev-type  
identities were recently studied in \cite{HanKW}, \cite{DelanoeKW}.

It is an interesting problem to characterize
the functions $f(x)$ which may arise as $\sigma_k$-curvature
functions in the conformal class of the round sphere.
An announcement of some work by Chang-Han-Yang in this direction for
$\sigma_2$ in dimension four was made in \cite{Hanlocal}.
For $k=1$, this problem is quite famous and has been studied in great depth. 
We do not attempt to make a complete list of references 
for this problem, we mention only \cite{CGYOld}, \cite{CSLIN}, 
\cite{EscobarSchoen}, \cite{Liscalar1}, \cite{Liscalar2}. 

\section{Admissible metrics}

 A natural question is: when does a manifold admit 
globally a strictly $k$-admissible metric, that is, 
a metric $g$ with $A_g \in \Gamma_k^+$?
In Section \ref{d4}, we have already discussed the beautiful 
result for $\sigma_2$ in dimension four by Chang-Gursky-Yang.
In the locally conformally flat case,
there are various topological restrictions. 
Guan-Lin-Wang showed 
\begin{theorem}(Guan-Lin-Wang \cite{GuanLinWang3})
Let $(M^n, g)$ be a compact, locally conformally flat 
manifold with $\sigma_1(A) > 0$. 

(i) If $A \in \overline{\Gamma}_k^+$ for some 
$2 \leq k < n/2$, then the $q$th Betti number $b_q = 0$
for 
\begin{align}
\left[ \frac{n+1}{2} \right] + 1 -k \leq q \leq
n - \left( \left[ \frac{n+1}{2} \right] +1 -k \right).
\end{align}
(ii) Suppose $ A \in \Gamma_2^+$, then $b_q = 0$ for
\begin{align}
\left[ \frac{n - \sqrt{n}}{2} \right] \leq q \leq 
\left[  \frac{n + \sqrt{n}}{2} \right].
\end{align}
 If $A \in \overline{\Gamma}_2^+$, $p =  \frac{n - \sqrt{n}}{2}$, 
and $b_p \neq 0 $, then $(M,g)$ is a quotient of 
$S^{n-p} \times H^p$. 

(iii) If $k \geq  \frac{n - \sqrt{n}}{2}$ and 
$A \in \Gamma_k^+$, then $b_q = 0$ for any
$2 \leq q \leq n-2$. If $k = \frac{n - \sqrt{n}}{2}$,
$A \in \overline{\Gamma}_k^+$, and $b_2 \neq 0$, then 
$(M,g)$ is a quotient of $S^{n-2} \times H^2$.
\end{theorem}

The proof of this theorem involves a careful analysis 
of the curvature terms in the Weitzenb\"ock formula 
for the Hodge Laplacian. For $\sigma_2$, the following was 
shown 

\begin{theorem}(Chang-Hang-Yang \cite{CHY})
Let $(M^n, g), (n \geq 5)$ be a smooth locally 
conformally flat Riemannian manifold such 
that $\sigma_1(A) > 0$, and $\sigma_2(A) \geq 0$, 
then $\pi_1(M) = 0$ for any $2 \leq i \leq [n/2] +1$, 
and $H^j(M, \mathbb{R}) = 0$ for any 
$n/2 - 1 \leq j \leq n/2 +1$. 
\end{theorem}
 Briefly, the technique in \cite{CHY} is to use the positivity 
condition to estimate the Hausdorf dimension of 
the singular set, using the developing map, as in 
\cite{SchoenYau}.
 Subsequently, in her thesis, Mar\'ia del Mar Gonz\'alez generalized 
this to prove the following. 
\begin{theorem}(Gonz\'alez \cite{Gonzalez1}) Let
$(M,g)$ be compact, locally conformally flat, 
and $A \in \Gamma_k^+$, $k < n/2$. 
Then for any $ 2 \leq i \leq \left[ \frac{n}{2} \right] + k -1$, 
the homotopy group $\pi_i(M) = \{0 \}$,
and the cohomology group $H^i(M, \mathbb{R}) = \{ 0 \}$ for 
$ \frac{n-2k}{2} + 1 \leq i \leq \frac{ n+2k}{2} -1$. 
\end{theorem}

 A beautiful gluing theorem was proved in \cite{GuanLinWang3}:
\begin{theorem}(Guan-Lin-Wang \cite{GuanLinWang3})
Let $2 \leq k < n/2$, and let $M_1^n$ and $M_2^n$ 
be two compact manifolds with $A_1, A_2 \in \Gamma_k^+$. 
Then the connected sum $M_1 \# M_2$ also admits a metric $g_{\#}$
with $A \in \Gamma_k^+$. If in addition, $M_1$ and $M_2$ are
locally conformally flat, then $g_{\#}$ can also be taken 
to be locally conformally flat. 
\end{theorem}
 This result can be viewed as a generalization of the 
analogous result for positive scalar curvature \cite{GromovLawson}
\cite{SYpos}. This yields many new examples of manifolds 
admitting metrics with $A \in \Gamma_k^+$, we refer the 
reader to \cite{GuanLinWang3} for more details. 

\section{The negative cone}
All of the previous results are concerned with the
positive curvature case. The negative curvature 
case exhibits quite different behaviour.
In this case, a serious technical difficulty arises in attempting to 
derive {\em a priori} second derivative estimates on solutions \cite{Jeff2}. 
Consider instead the following generalization of the 
Schouten tensor. Let $t \in \mathbb{R}$, and define 
\begin{align}
A^t = \frac{1}{n-2} \left( Ric - \frac{t}{2(n-1)}Rg \right).
\end{align}
We let $\Gamma_k^- = - \Gamma_k^+$. 
\begin{theorem}(Gursky-Viaclovsky \cite{GVNegative})
\label{sharpthm}Let $(M,g)$ be a compact Riemannian manifold, 
assume that $A^t_g \in \Gamma_k^{-}$
for some $t<1$, and let $f(x) < 0$ be any smooth 
function on $M^n$. Then there 
exists a unique conformal metric $\tilde{g} = 
e^{2w}g$ satisfying
\begin{align}
\label{eqn1}
\sigma_k^{1/k}(\tilde{g}^{-1}A^t_{\tilde{g}}) = f(x).
\end{align}
\end{theorem}
As noted above, the 
second derivative estimate encounters technical difficulties for $t = 1$, 
but for $t<1$, this difficulty can be overcome.
Also, for $t = 1 + \epsilon$, the equation is not necessarily elliptic, 
therefore $t=1$ is critical for more than one reason.
Some local counterexamples to the 
second derivative estimate in the negative case have been
given in \cite{STW}, but there are
no known global counterexamples. 
The above theorem was subsequently proved by 
parabolic methods in \cite{LiSheng}. 

 The above theorem has the following corollary.
Using results of \cite{Brooks}, \cite{GaoYau}, and \cite{Lohkamp}, 
every compact manifold of dimension $n \geq 3$ admits a metric
with negative Ricci curvature. Therefore applying the 
theorem when $t=0$,
\begin{corollary}(Gursky-Viaclovsky \cite{GVNegative})
Every smooth compact $n$-manifold, $n \geq 3$, admits a Riemannian metric with
$Ric < 0$ and  
\begin{align}
\mbox{\em{det}}(g^{-1}Ric) = \mbox{\em{constant}}. 
\end{align}
\end{corollary} 

 It turns out the equation is 
also elliptic for $t \geq n-1$, and this has some 
interesting consequences, see \cite{ShengZhang} for details. 
We also mention that Mazzeo and Pacard considered the 
$\sigma_k$-Yamabe equation in the context of conformally compact 
metrics, and showed that the deformation problem 
is unobstructed \cite{MazzeoPacard}. 

\bibliography{Viaclovsky_references}
\end{document}